\input amstex
\documentstyle{amsppt}
\magnification=\magstep1
 \hsize 13cm \vsize 18.35cm \pageno=1
\loadbold \loadmsam
    \loadmsbm
    \UseAMSsymbols
\topmatter
\NoRunningHeads
\title On the Euler Numbers and Its Applications
\endtitle
\author
  Taekyun Kim
\endauthor

\abstract Recently, the $q$-Euler numbers and polynomials are
constructed in [T. Kim, The modified $q$-Euler numbers and
polynomials, Advanced Studies in Contemporary Mathematics, 16(2008),
161-170]. These $q$-Euler numbers and polynomials contain the
interesting properties. In this paper we prove Von-Staudt Clausen's
type theorem related to the $q$-Euler numbers. That is, we prove
that the $q$-Euler numbers are $p$-adic integers. Finally, we give
the proof of Kummer type congruences for the $q$-Euler numbers.
\endabstract
\thanks  2000 AMS Subject Classification: 11B68, 11S80
\endthanks
\endtopmatter

\document

{\bf\centerline {\S 1. Introduction/Definition}}

 \vskip 20pt
Let $p$ be a fixed odd prime. Throughout this paper $\Bbb Z_p ,$
$\Bbb Q_p ,$ $\Bbb C,$ and $\Bbb C_p$ will, respectively, denote the
ring of $p$-adic rational integers, the field of $p$-adic rational
numbers, the complex number field, and the completion of algebraic
closure of $\Bbb Q_p .$  Let $v_p$ be the normalized  exponential
valuation of $\Bbb C_p$ with $|p|_p=p^{-v_p(p)}=\frac{1}{p}.$ When
one talks of $q$-extension, $q$ is variously considered as an
indeterminate, a complex $q \in \Bbb C$, or a $p$-adic number $q\in
\Bbb C_p$, see [9-22]. If $q\in \Bbb C$, then we assume $|q|<1.$ If
$q\in \Bbb C_p$, then we assume $|1-q|_p<p^{-\frac{1}{p-1}}.$ The
ordinary Euler numbers are defined as
$$\frac{2}{e^t +1}=e^{Et}=\sum_{n=0}^{\infty}E_n\frac{t^n}{n!}, $$
where we use the technique method notation by replacing $E^n$ by
$E^n$ ($n\geq 0$), symbolically (see [1-23]). From this definition,
we can derive the following relation.
$$E_0=1, \text{ and  } (E+1)^n+E_n=2\delta_{0,n}, \text{
where $\delta_{0,n}$ is Kronecker symbol.} $$ For $x\in\Bbb Q_p $ (
or $\Bbb R$), we use the notation $[x]_q=\frac{1-q^x}{1-q}, $ and
$[x]_{-q}=\frac{1-(-q)^x }{1+q},$ see [5-6]. In [5], the $q$-Euler
numbers are defined as
$$E_{0,q}=\frac{[2]_q}{2}, \text{ and }
(qE+1)^n+E_{n,q}=[2]_q\delta_{0,n},  \tag1$$ with the usual
convention of replacing $E^n$ by $E_{n,q}.$
Note that
$\lim_{q\rightarrow 1}E_{n,q}=E_n .$ For a fixed positive integer
$d$ with $(p,d)=1$, let
$$X=X_d=\varprojlim_N \Bbb Z/dp^N\Bbb Z , \;\;X_1=\Bbb Z_p,$$
$$X^*=\bigcup\Sb 0<a<dp\\ (a,p)=1\endSb a+dp\Bbb Z_p,$$
$$a+dp^N\Bbb Z_p=\{x\in X\mid x\equiv a\pmod{dp^N}\},$$
where $a\in \Bbb Z$ lies in $0\leq a<dp^N$, (see [4-23]). We say
that $f$ is a uniformly differentiable function at a point $a
\in\Bbb Z_p $ and denote this property by $f\in UD(\Bbb Z_p )$, if
the difference quotients $F_f (x,y) = \dfrac{f(x) -f(y)}{x-y} $ have
a limit $l=f^\prime (a)$ as $(x,y) \to (a,a)$.
 For $f\in UD(\Bbb Z_p
)$, let us start with the expression
$$\eqalignno{ & \dfrac{1}{[p^N ]_q} \sum_{0\leq j < p^N} q^j f(j)
=\sum_{0\leq j < p^N} f(j)
\mu_q (j +p^N \Bbb Z_p ), }$$ representing a $q$-analogue of Riemann
sums for $f$, see [5, 6]. The integral of $f$ on $\Bbb Z_p$ will be
defined as limit ($n \to \infty$) of those sums, when it exists. The
$q$-deformed bosonic $p$-adic integral of the function $f\in UD(\Bbb
Z_p )$ is defined by
$$ I_q (f)= \int_{\Bbb Z_p }f(x) d\mu_q (x) = \lim_{N\to \infty}
\dfrac{1}{[dp^N ]_q} \sum_{0\leq x < dp^N} f(x) q^x, \text{ see [5]}
. \tag2 $$ In the sense of fermionic, let us define the fermionic
$p$-adic $q$-integral as
$$I_{-q}(f)=\int_{\Bbb Z_p}f(x)d\mu_{-q}(x)=\lim_{N\rightarrow
\infty}\frac{1}{[p^N]_{-q}}\sum_{x=0}^{p^N-1}f(x)(-q)^x, \text{ see
[5-10]} .\tag3$$ From (3) we note that
$$qI_{-q}(f_1)+I_{-q}(f)=[2]_qf(0), \text{ where $f_1(x)=f(x+1)$.}
\tag4$$ In [5], the Witt's type formula  for the $q$-Euler numbers
$E_{n,q}$ id given by $$\sum_{n=0}^{\infty}\int_{\Bbb
Z_p}[x]_q^nq^{-x}d\mu_{-q}(x)\frac{t^n}{n!}=\int_{\Bbb Z_p}
e^{[x]_qt}q^{-x}d\mu_{-q}(x)=\sum_{n=0}^{\infty}E_{n,q}\frac{t^n}{n!}.
\tag5$$ By comparing the coefficients on both sides  in (5), we see
that
$$\int_{\Bbb Z_p}[x]_q^n q^{-x} d\mu_{-q}(x)=E_{n,q}, \text{ see
[5]}. \tag6$$ By the definition of the fermionic $p$-adic
$q$-integral on $\Bbb Z_p$, the $q$-Euler polynomials are also
defined as
$$\int_{\Bbb Z_p}e^{[x+y]_qt}q^{-y}d\mu_{-q}(y)=e^{[x]_qt}\int_{\Bbb
Z_p}e^{q^x[y]_qt}q^{-y}d\mu_{-q}(y)=\sum_{n=0}^{\infty}E_{n,q}(x)\frac{t^n}{n!}.
\tag7$$ From (6) and (7), we note that
$$E_{n,q}(x)=\sum_{k=0}^n \binom{n}{k}[x]_q^{n-k}q^kE_{k,q}, \text{ where $\binom{n}{k}=
\frac{n (n-1)\cdots (n-k+1)}{k!}$}.\tag8$$ Let $F_q(t,x)$ be the
generating function of the $q$-Euler polynomials. Then we have
$$F_q(t,x)=\sum_{n=0}^{\infty}E_{n,q}(x)\frac{t^n}{n!}=[2]_q\sum_{k=0}^{\infty}(-1)^ke^{[k+x]_qt}.
\tag9$$ Let $\chi$ be the Dirichlet's character with odd conductor
$d\in\Bbb N$. Then the generalized $q$-Euler numbers attached to
$\chi$ are defined as
$$E_{n,\chi,q}=\int_{X}[x]_q^n q^{-x} \chi(x)d\mu_{-q}(x)=[d]_q^n
\frac{[2]_q}{[2]_{q^d}}\sum_{a=0}^{d-1}\chi(a)(-1)^aE_{n,q^d}(\frac{a}{d}).
\tag10$$ Let
$F_{\chi,q}(t)=\sum_{n=0}^{\infty}E_{n,\chi,q}\frac{t^n}{n!}.$ Then
we note that
$$F_{\chi,q}(t)=\sum_{n=0}^{\infty}E_{n,\chi,q}\frac{t^n}{n!}=[2]_q\sum_{n=0}^{\infty}\chi(n)(-1)^ne^{[n]_qt}.\tag11$$
In this paper we prove the Von-Staudt-Clausen's  type theorem
related to the $q$-Euler numbers. That is, we prove that the
$q$-Euler numbers are the $p$-adic integers. Finally, we give the
proofs of the  Kummer congruences for the $q$-Euler numbers.

\vskip 10pt

{\bf\centerline {\S 2. $q$-Euler numbers and polynomials}}
 \vskip 10pt
 From (1)and (6) we derive
$$E_{n,q}=\int_{\Bbb Z_p}q^{-x}[x]_q^n
d\mu_{-q}(x)=\frac{[2]_q}{2}\int_{\Bbb Z_p}[x]_q^nd\mu_{-1}(x).$$
Thus, we note that  $\lim_{n\rightarrow
\infty}E_{n,q}=E_n=\int_{\Bbb Z_p}x^n d\mu_{-1}.$ For $q\in\Bbb C_p$
with $|1-q|_p<p^{- \frac{1}{p-1}} ,$ we have
$$(-1)^j[j]_q-j(-1)^j=(-1)^j \left(\frac{\sum_{l=0}^j \binom{j}{l}(q-1)^l -1}{q-1}-j\right)
=(-1)^j\sum_{l=2}^j\binom{j}{l}(q-1)^{l-1} .$$ Thus, we see that
$$\left|(-1)^j([j]_q-j)\right|_p \leq \max_{2\leq l \leq j} \left( |(q-1)|_p^{l-1}\right)=|q-1|_p. \tag12 $$
From (12), we can derive
$$\aligned \left|\sum_{j=0}^{p-1}(-1)^j[j]_q\right|_p&=
\left|\sum_{j=0}^{p-1}(-1)^j([j]_q-j)+\sum_{j=0}^{p-1}(-1)^jj\right|_p\\
&=\left|\sum_{j=0}^{p-1}(-1)^j([j]_q-j)+\frac{p-1}{2}\right|_p \leq
1.\endaligned\tag13$$
For $k\geq 1$, let
 $$T_n(k)=\sum_{x=0}^{p^k-1}(-1)^x[x]_q^n=[0]_q^n-[1]_q^n+\cdots+[p^k-1]_q^n.
 \tag14$$
 Note that $\lim_{k\rightarrow
 \infty}T_n(k)=\frac{2}{[2]_q}E_{n,q}.$
 From (14), we can derive
 $$\aligned
&T_{n}(k+1)=\sum_{x=0}^{p^{k+1}-1}(-1)^x[x]_q^n=\sum_{i=0}^{p^k-1}\sum_{j=0}^{p-1}
[i+jp^k]_q^n(-1)^{i+jp^k}\\
&=\sum_{i=0}^{p^k-1}\sum_{j=0}^{p-1}\left([i]_q+q^i[jp^k]_q\right)^n(-1)^{i+jp^k}
=\sum_{i=0}^{p^k-1}\sum_{j=0}^{p-1}\sum_{l=0}^n\binom{n}{l}[i]_q^{n-l}q^{il}[jp^k]_q^l(-1)^{i+jp^k}\\
&=\sum_{i=0}^{p^k-1}\sum_{j=0}^{p-1}\sum_{l=0}^n\binom{n}{l}[i]_q^{n-l}q^{il}[p^k]_q^l
[j]_{q^{p^k}}^l(-1)^{i+j}\\
&=\sum_{i=0}^{p^k-1}[i]_q^n(-1)^i+\sum_{i=0}^{p^k-1}\sum_{j=0}^{p-1}\sum_{l=1}^n\binom{n}{l}[i]_q^{n-l}
q^{il}[p^k]_q^l[j]_{q^{p^k}}^l.
 \endaligned\tag15$$
 Thus, we have
 $$T_n(k+1)\equiv \sum_{i=0}^{p^k-1}[i]_q^n(-1)^i  \text{ ($\mod [p^k]_q$ )}. \tag16$$
 From (15) we note that
 $$\aligned
\sum_{x=0}^{p^{k+1}-1}[x]_q^n(-1)^x&=\sum_{j=0}^{p-1}\sum_{i=0}^{p^k-1}[j+ip]_q^n(-1)^{j+ip}\\
&=\sum_{j=0}^{p-1}(-1)^j\sum_{i=0}^{p^k-1}(-1)^i\sum_{l=0}^n\binom{n}{l}[j]_q^{n-l}q^{jl}[p]_q^l[i]_{q^p}^l\\
&=\sum_{j=0}^{p-1}(-1)^j[j]_q^n+\sum_{j=0}^{p-1}(-1)^j\sum_{i=0}^{p^k-1}(-1)^i\sum_{l=1}^n\binom{n}{l}
[j]_q^{n-l}q^{jl}[p]_q^l[i]_{q^p}^l\\
&\equiv \sum_{j=0}^{p-1}(-1)^{j}[j]_q^n \text{ ($\mod [p]_q $ )}.
 \endaligned\tag17$$
By (17), we obtain
$$\sum_{x=0}^{p^{k+1}-1}(-1)^x[x]_q^n\equiv
\sum_{j=0}^{p-1}(-1)^j[j]_q^n \text{ ($\mod [p]_q$ ) }. \tag18$$
From (16) and (19), we can also derive
$$\sum_{j=0}^{p-1}(-1)^j[j]_q^n\equiv
\frac{2}{[2]_q}\int_{X}[x]_q^nq^{-x}d\mu_{-q}(x)=\frac{2}{[2]_q}E_{n,q}
\text{ ($\mod [p]_q  $) }. \tag19$$ Thus, we note that
$$\sum_{j=0}^{p-1}(-1)^j[j]_q^n\equiv
\frac{2}{[2]_q}E_{n,q} \text{ ($\mod [p]_q $) }. \tag20$$ Therefore
we obtain the following theorem.

\proclaim{ Theorem 1} For $n\geq 0$, we have
$$\sum_{j=0}^{p-1}(-1)^j[j]_q^n\equiv
\frac{2}{[2]_q}E_{n,q} \text{ $(\mod [p]_q )$ }. $$
\endproclaim
By (15), (16) and (20), we obtain the following corollary.

\proclaim{ Corollary 2} For $n\geq 0$, we have
$$\frac{2}{[2]_q}E_{n,q}+\sum_{j=0}^{p-1}(-1)^{j+1}[j]_q^n \in \Bbb
Z_p.$$
\endproclaim
For $n\geq 0$, we note that
$$\aligned
&\left|\frac{2}{[2]_q}E_{n,q}\right|_p=\left|\frac{2}{[2]_q}E_{n,q}-\sum_{j=0}^{p-1}(-1)^j[j]_q^n+\sum_{j=1}^{p-1}(-1)^j[j]_q^n\right|_p\\
& \leq \max\left(
\left|\frac{2}{[2]_q}E_{n,q}-\sum_{j=0}^{p-1}(-1)^j[j]_q^n\right|_p,
\text{ }
  \left|\sum_{j=1}^{p-1}(-1)^j[j]_q^n\right|_p \right).\endaligned$$
By (13) and Corollary 2, we obtain the following corollary.

\proclaim{ Corollary 3} For $n\geq 0$, we have
$$\frac{2}{[2]_q}E_{n,q} \in \Bbb
Z_p.$$
\endproclaim

Let $\chi$ be the Dirichlet's character with odd conductor $d
(\in\Bbb N)$. Then the generalized $q$-Euler numbers attached to
$\chi$ as follows.
$$\sum_{n=0}^{\infty}E_{n,\chi,q}\frac{t^n}{n!}=[2]_q\sum_{n=0}^{\infty}\chi(n)(-1)^ne^{[n]_qt}
=\int_{X}\chi(x)e^{[x]_qt}q^{-x}d\mu_{-q}(x).\tag21$$
 We denote $\bar{d}=[d, p]$ the least common multiple of the conductor
 $d$ of $\chi$ and $p$. From (21), we derive
 $$\frac{2}{[2]_q}E_{n,\chi,q}=\frac{2}{[2]_q}\int_{X}[x]_q^nq^{-x}\chi(x)
 d\mu_{-q}(x)=\lim_{N\rightarrow
 \infty}\sum_{x=0}^{dp^N-1}[x]_q^n\chi(x)(-1)^x. \tag22 $$
By (22), we see that
 $$\aligned
&\frac{2}{[2]_q}E_{n,\chi,q}=\lim_{\rho \rightarrow \infty}\sum\Sb
1\leq x \leq\bar{d}p^{\rho}\\ (x,p)=1\endSb \chi(x)(-1)^x[x]_q^n+
\lim_{\rho\rightarrow \infty}\sum_{y=1}^{\bar{d}p^{\rho-1}}\chi(p)\chi(y)[p]_q^n[y]_{q^p}^n(-1)^y\\
&=\lim_{\rho \rightarrow \infty}\sum\Sb 1\leq x
\leq\bar{d}p^{\rho}\\ (x,p)=1\endSb
\chi(x)(-1)^x[x]_q^n+\chi(p)[p]_q^n \lim_{\rho\rightarrow
\infty}\sum_{y=1}^{\bar{d}p^{\rho-1}}\chi(y)[y]_{q^p}^n(-1)^y .
 \endaligned$$
 Thus, we have
$$\frac{2}{[2]_q}E_{n,\chi,q}-\chi(p)[p]_q^n\frac{2}{[2]_{q^p}}E_{n,\chi,
q^{p}}= \lim_{\rho \rightarrow \infty}\sum\Sb 1\leq x
\leq\bar{d}p^{\rho}\\ (x,p)=1\endSb \chi(x)(-1)^x[x]_q^n .\tag23$$
Let $w$ denote the Teichm$\ddot{u}$ller character $\mod p$. For
$x\in X^{*}$, we set $<x>=<x:q>=\frac{[x]_q}{w(x)}.$ Note that
$|<x>-1|_p<p^{-\frac{1}{p-1}},$ where $<x>^s=\exp(s\log_p <x> )$ for
$s\in\Bbb Z_p .$ For $s \in \Bbb Z_p$,  we define the $p$-adic
$q$-$L$-function related to $E_{n,\chi,q}$ as follows.
$$L_{p,q,E}(s, \chi)=\lim_{\rho\rightarrow \infty}\sum\Sb
1\leq x \leq\bar{d}p^{\rho}\\ (x,p)=1\endSb \chi(x)(-1)^x<x>^{-s}.
\tag24$$
For $k \geq 0 ,$ we have
$$\aligned
&L_{p,q,E}(-k, \chi w^k)=\lim_{\rho\rightarrow \infty}\sum\Sb 1\leq
x \leq\bar{d}p^{\rho}\\ (x,p)=1\endSb \chi(x)(-1)^x [x]_q^k \\
&=\frac{2}{[2]_q}\int_{X}[x]_q^k\chi(x)q^{-x} d\mu_{-q}(x)
-\frac{2}{[2]_{q^p}}\int_{pX}[x]_q^k\chi(x)q^{-x} d\mu_{-q}(x)\\
&=\frac{2}{[2]_q}\int_{X}[x]_q^k\chi(x)q^{-x} d\mu_{-q}(x)
-\chi(p)[p]_q^k\frac{2}{[2]_{q^p}}\int_{X}[x]_{q^p}^k\chi(x)q^{-px} d\mu_{-q^p}(x)\\
&=\frac{2}{[2]_q}E_{n,\chi,q}-\chi(p)[p]_q^k\frac{2}{[2]_{q^p}}E_{n,\chi,
q^p}. \endaligned\tag25$$

It is easy to see that $<x>^{p^n}\equiv 1 $ ($\mod p^n$ ). From the
definition of $L_{p,q, E}(s, \chi) ,$ we can derive
$$\aligned
&L_{p,q,E}(-k, \chi )=\lim_{\rho\rightarrow \infty}\sum\Sb 1\leq x \leq\bar{d}p^{\rho}\\
(x,p)=1\endSb \chi(x)(-1)^x <x>^k \\
&\equiv \lim_{\rho\rightarrow \infty}\sum\Sb 1\leq x \leq\bar{d}p^{\rho}\\
(x,p)=1\endSb \chi(x)(-1)^x <x>^{k^{\prime}} \text{ ( $\mod p^n $
)},
\endaligned$$
whenever $ k \equiv k^{\prime} $ ($\mod p^n(p-1)$ ). That is,
$L_{p,q,E}(-k, \chi w^k) \equiv L_{p,q,E}(-k^{\prime}, \chi
w^{\prime})$ ($\mod p^n $).

Therefore we obtain the following theorem.

\proclaim{Theorem 4} (Kummer Congruence ) For $k \equiv k^{\prime}$
$(\mod p^n(p-1))$ , we have
 $$ \frac{2}{[2]_q}E_{k,\chi,q}-\frac{2}{[2]_{q^p}}E_{k,\chi,
 q^p}\equiv \frac{2}{[2]_q}E_{k^{\prime},
 \chi,q}-\frac{2}{[2]_{q^p}} E_{k^{\prime},\chi,q^p} \text{ $(\mod p^n  )$}. $$
\endproclaim

Let $\chi$ be the primitive Dirichlet's character with conductor
$p$. Then we have
$$\aligned
&\sum_{x=0}^{p^{N+1}-1}\chi(x)(-1)^x[x]_q^n=
\sum_{a=0}^{p-1}\sum_{x=0}^{p^N-1} \chi(a+px)(-1)^{a+px}[a+px]_q^n \\
&=\sum_{a=0}^{p-1}\chi(a)(-1)^a\sum_{x=0}^{p^N-1}(-1)^x([a]_q+q^a[p]_q[x]_{q^p})^n\\
&=\sum_{a=0}^{p-1}\chi(a) (-1)^a
\sum_{x=0}^{p^N-1}(-1)^x\sum_{l=0}^n
\binom{n}{l}[a]_q^{n-l}q^{al}[p]_q^l[x]_{q^p}^l\\
&\equiv \sum_{a=0}^{p-1}\chi(a)(-1)^a[a]_q^n \text{ ($\mod [p]_q$
)}. \endaligned$$
If $\rho \rightarrow \infty$, then we have
$$\frac{2}{[2]_q}\int_{X}\chi(x)(-1)^x[x]_q^nq^{-x}d\mu_{-q}(x)\equiv
\sum_{a=0}^{p-1}\chi(a)(-1)^a[a]_q^n  \text{ ($\mod [p]_q $ )}.$$
Thus, we can obtain the following. Let $\chi$ be the primitive
Dirichlet's character with conductor $p$. Then we have
$$\frac{2}{[2]_q}E_{n,\chi,q}\equiv
\sum_{a=0}^{p-1}\chi(a)(-1)^a[a]_q^n \text{ ($\mod [p]_q$ ) }.
\tag26$$

The Eq.(26) also seems to be the new interesting formula.  As
$q\rightarrow 1$, we can also obtain
$$ E_{n,\chi}\equiv \sum_{a=0}^{p-1}\chi(a)(-1)^a a^n \text{ ($\mod
p $ )}.$$

\vskip 20pt

\quad Taekyun Kim

\quad Division of General Education-Mathematics, Kwangwoon
University, Seoul 139-701, S. Korea \quad e-mail:\text{
tkkim$\@$kw.ac.kr}

\enddocument